\documentclass[a4paper,12pt]{amsart}
\usepackage{amsmath}
\usepackage{amssymb}
\usepackage{mathrsfs}
\usepackage{enumerate}
\usepackage{ifthen}
\usepackage{graphicx}
\usepackage{caption}
\usepackage[T1]{fontenc} 
\usepackage{tikz}
\usepackage{cite}

\nonstopmode \numberwithin{equation}{section}
\newtheorem{theorem}{Theorem}[section]

\newtheorem{corollary}{Corollary}[section]

\newtheorem{lemma}{Lemma}[section]

\theoremstyle{remark}
\theoremstyle{definition}
\newtheorem{remark}{Remark}[section]
\newtheorem{definition}{Definition}[section]

\theoremstyle{plain}

\newtheorem*{lemA}{Lemma A}
\newtheorem*{lemB}{Lemma B}
\newtheorem*{lemC}{Lemma C}

\numberwithin{equation}{section}
\numberwithin{theorem}{section}

\newcounter{minutes}
\divide\time by 60
\newcounter{hours}
\multiply\time by 60
\addtocounter{minutes}{-\time}

\begin{document}

\title{Radius of concavity for certain classes of functions}

\author{Molla Basir Ahamed$^*$}
\address{Molla Basir Ahamed, Department of Mathematics, Jadavpur University, Kolkata-700032, West Bengal, India.}
\email{mbahamed.math@jadavpuruniversity.in}

\author{Rajesh Hossain}
\address{Rajesh Hossain, Department of Mathematics, Jadavpur University, Kolkata-700032, West Bengal, India.}
\email{rajesh1998hossain@gmail.com}



\begin{abstract}
Let $\mathcal{S}(p)$ be the class of all meromorphic univalent functions defined in the unit disc $\mathbb{D}$ of the complex plane with a simple pole at $z=p$ and normalized by the conditions $f(0)=0$ and $f^{\prime}(0)=1$. In this paper, we find radius of concavity and compute the same for functions in $\mathcal{S}(p)$ and for some other well-known classes of functions on unit disk. We explore general linear combinations $F(z):=\lambda_{1}f_{1}(z)+\cdot\cdot\cdot+\lambda_{2n}f_{2n}(z)$, $\lambda_{j}\in\mathbb{C}$, $n\in\mathbb{N}$ of functions belonging to the class $\mathcal{S}(p)$ and some other classes of functions of analytic univalent functions and investigate their radii of univalence, convexity and concavity.
\end{abstract}
\maketitle
\pagestyle{myheadings}
\markboth{M. B. Ahamed and R. Hossain}{Radius of concavity for certain class of functions}
\section{\bf Introduction}
Throughout the article, we will denote the open unit disk of the complex plane $ \mathbb{C} $ by $ \mathbb{D}:=\{z\in\mathbb{C} : |z|<1\} $. Let $\mathcal{A}$ be the class of functions $f$ analytic in the unit disk $\mathbb{D}$ and normalized by $f(0)=0$ and $f^{\prime}(0)=1$. We denote by $\mathcal{S}$ the class of univalent functions in $\mathcal{A}$. For $\alpha\leq 1$, let $S^*(\alpha)$ denote the subclass of $\mathcal{A}$ of starlike functions of order $\alpha$. If $0\leq \alpha<1$, then $f$ is said to be univalent and starlike of order $\alpha$, \textit{i.e.,} we have $\mathcal{S}^*_{\alpha}\subset \mathcal{S}$. A function $f\in\mathcal{A}$ is said to be starlike of order $\alpha$ if, and only if, 
\begin{align*}
	{\rm Re} \left(\frac{z f^{\prime}(z)}{f(z)}\right)>\alpha.
\end{align*}
The class of starlike functions of order $\alpha$ is denoted by $S^*(\alpha)$.  \vspace{1.2mm}

Let $\mathcal{C}_{\alpha}$ be the class of convex functions of order $\alpha$ corresponding to $\mathcal{S}_{\alpha}$ in the usual way that $g\in\mathcal{C}_{\alpha} \Leftrightarrow zg^{\prime}\in \mathcal{S}^*(\alpha)$. Note that $\mathcal{S}^*(0)=\mathcal{S}^*$ and $\mathcal{C}_{0}=\mathcal{C}$. It is well-known that $f\in\mathcal{C}$ if, and only if
\begin{align*}
	{\rm Re}\left(1+\frac{z f^{\prime\prime}(z)}{f^{\prime}(z)}\right)>0,\; z\in\mathbb{D},
\end{align*}
and $f\in\mathcal{S}^*$ if, and only if
\begin{align*}
	{\rm Re}\left(\frac{z f^{\prime}(z)}{f(z)}\right)>0,\; z\in\mathbb{D}.
\end{align*}
Further, it is worth mentioning that $\mathcal{C}\subsetneq\mathcal{S}^*$. Interestingly, each $f\in\mathcal{S}$ maps $\mathbb{D}_r:=\{z\in\mathbb{C} : |z|<r\}$, $r\in (0, 1]$ onto a convex domain if, and only if, ${\rm Re}\left(1+zf^{\prime\prime}(z)/f^{\prime}(z)\right)>0$ for $z\in\mathbb{D}_r$ (see \cite[p. 105, Problem 108]{Pólya-Szeg˝o-1954}).\vspace{1.2mm}

We now turn to defining the class of functions that is our main focus in this article. 
\begin{definition}
	A function $f$ is said to belong to the class of \emph{concave functions with an opening angle at infinity less than or equal to $\pi A$, $A\in [1, 2]$} if it satisfies the following conditions:
	\begin{enumerate}
		\item[(i)] $f$ is analytic and univalent in $\mathbb{D}$.
		\item[(ii)] $f$ maps $\mathbb{D}$ conformally onto a set whose complement with respect to $\mathbb{C}$ is convex and satisfies the normalization $f(0)=0=f^{\prime}(0)-1$ and $f(1)=\infty$.
		\item[(iii)] The opening angle of $f(\mathbb{D})$ at infinity is less than or equal to $\pi A$, $A\in [1,2].$
	\end{enumerate}
\end{definition}
We denote the family of concave univalent functions with opening angle $\pi A$, $A\in [1,2],$ at infinity by $\widehat{{\rm Co}}(A)$. It is important to note that the boundary of $f(\mathbb{D})$, $f\in \widehat{{\rm Co}}(A)$, is contained in a wedge shaped region with opening angle $\pi A$ but not in any bigger opening angle. However, it should be observed that for $f\in \widehat{{\rm Co}}(A)$, $A\in [1, 2]$, the closed set $\mathbb{C}\setminus f(\mathbb{D})$ is convex and unbounded. Furthermore, it is clear that $\widehat{{\rm Co}}(2)$ contains the classes $\widehat{{\rm Co}}(A)$, $A\in [1, 2]$ (see \cite{Avkhadiev-Wirths-JM-2005}). Furthermore, the image domain of functions belonging to the class \( \widehat{{\rm Co}}(1) \) becomes a convex half-plane, hence, we can see that ${\rm Co}(A)=\widehat{{\rm Co}}(A)\cap\mathcal{S}$. In $2005$, Avkhadiev and Wirths \cite{Avkhadiev-Wirths-JM-2005} provided an analytic characterization of functions in the class \( {\rm Co}(A) \). The Fekete-Szeg\"o problem for this class was solved by Bhowmik \emph{et al.} \cite{Bhowmik-Ponnusamy-Wirths-JMAP-2011}. Abu-Muhanna and Ponnusamy \cite{Abu-Muhanna-Ponnusamy-MN-2017} discussed results related to the Yamashita conjecture on Dirichlet finite integrals for functions in \( {\rm Co}(A) \). For a detailed discussion on concave univalent functions, see \cite{Avkhadiev-CMH-2006,Bhowmik-MN-2012,Bhowmik-Das-JMAA-2018,Cruz-Pommerenke-CVEE-2007} and references therein.\vspace{2mm}

\begin{figure}
	\begin{center}
		\includegraphics[width=0.55\linewidth]{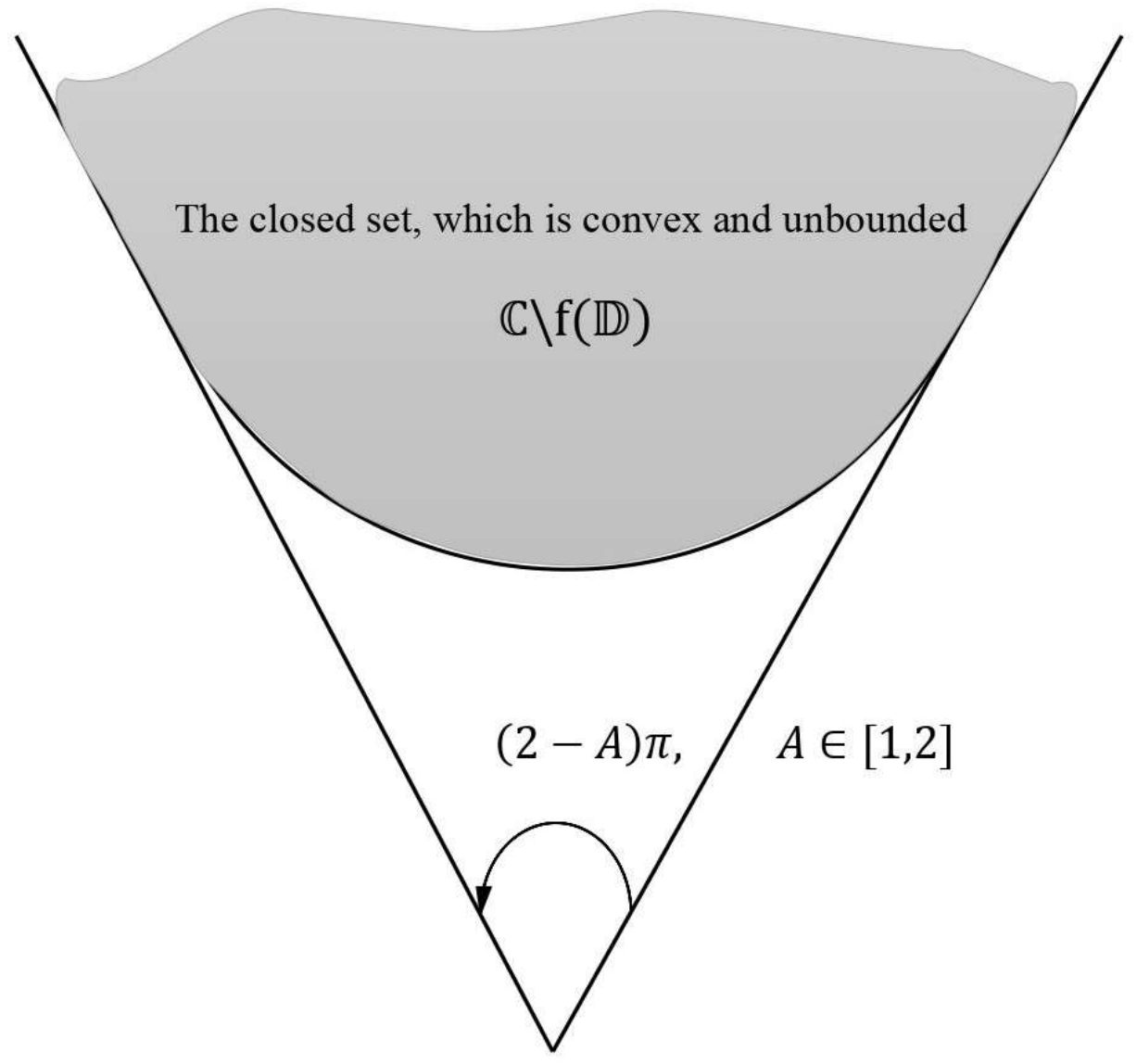}
	\end{center}
	\caption{The figure depicts the image domain of the concave function $ f $, characterized by an open angle of $ \pi A $ for $ A \in [1,2] $. The shaded region $ \mathbb{C}\setminus f(\mathbb{D}) $ denotes a closed set that is both convex and unbounded.}
\end{figure}

We consider the class $ \mathcal{A}(p) $, where $ p\in (0, 1) $ consisting of all meromorphic functions in $ \mathbb{D} $ with a simple pole at $ z=p $ and normalized by the condition $ f(0)=0=f^{\prime}(0)-1$. Let $ \mathcal{S}(p) $ be the class of all univalent functions in $ \mathcal{A}(p) $. In $ 1931 $, Fenchel \cite{Fenchel-1931} obtained the sharp lower bound of $ |f(z)| $ for $ f\in\mathcal{S}(p) $ and the upper bound was established by Kirwan and Schober (see \cite{Kirwan-Schober-JAM-1976}) in $ 1976 $. It is well-known that if $ f\in\mathcal{S}(p) $, then 
\begin{align}\label{Eq-1.1}
	|k_p(-r)|\leq |f(z)|\leq |k_p(z)|,\;\; |z|=r<1,
\end{align}
where 
\begin{align*}
	k_p(z):=\frac{-pz}{(z-p)(1-pz)}\in \mathcal{S}(p)\; \mbox{for}\; z\in\mathbb{D}.
\end{align*}
Let $ {\rm Co}(p) $ be a subclass of $ \mathcal{S}(p) $ consisting of functions $ f $ such that $ \overline{\mathbb{C}}\setminus f(\mathbb{D}) $ is a bounded convex set, where $ \overline{\mathbb{C}}:=\mathbb{C}\cup\{\infty\} $. In \cite{Pfaltzgraff-Pinchuk-JAM-1971}, it is proved that $ f\in {\rm Co}(p) $ if, and only if, $ f\in\mathcal{S}(p) $ and the function $ P_f : \mathbb{D}\to\mathbb{C} $ defined as
\begin{align}\label{Eq-1.2}
	P_f(z)=-\left(1+\frac{zf^{\prime\prime}(z)}{f^{\prime}(z)}+\frac{z+p}{z-p}-\frac{1+pz}{1-pz}\right),\; z\in\mathbb{D}\setminus\{p\}
\end{align}
with
\begin{align*}
	P_f(p):=\frac{1+p^2}{1-p^2},
\end{align*}
is analytic and ${\rm Re}\; P_f(z)>0$ for $z\in\mathbb{D}$.\vspace{2mm}

The notion of the radius of concavity was recently introduced in \cite{Bhowmik-CMFT-2024}. In this context, it is worth observing that for a function $f \in \mathrm{Co}(p)$ and $r \in (0,1)$, the scaled mapping $g(z) := r^{-1}f(rz)$ ($z \in \mathbb{D}$) does not necessarily belong to $\mathrm{Co}(p)$. Consequently, the convexity of $\overline{\mathbb{C}} \setminus f(\mathbb{D})$ does not automatically imply the convexity of $\overline{\mathbb{C}} \setminus f(\mathbb{D}_r)$ for each $0 < r \le 1$. This fact motivates the definition of the radius of concavity (with respect to $\mathrm{Co}(p)$) for subsets of $\mathcal{A}(p)$ in the following manner.
\begin{definition}\cite{Bhowmik-CMFT-2024}\label{Def-1.1}
The radius of concavity (w.r.t $ {\rm Co}(p) $) of a subset $ \mathcal{A}_1(p) $ of $ \mathcal{A}(p) $ is the largest number $ R_{{\rm Co}(p)}\in (0,1] $  such that for each function $ f\in\mathcal{A}_1(p) $, $ {\rm Re}\; P_f(z)>0 $ for all $ |z|< R_{\rm Co(p)}$, where $ P_{f}(z) $ is defined in \eqref{Eq-1.2}.
\end{definition}
The main objective of this paper is to obtain a lower bound for the radius of concavity $R_{\mathrm{Co}(p)}$ of the class $\mathcal{S}(p)$. In addition, we study the family $\mathrm{Co}(A)$ consisting of functions $f \in \mathcal{A}$ that map $\mathbb{D}$ conformally onto a domain with a convex complement in $\mathbb{C}$, normalized by $f(1) = \infty$. It is well known that $f \in \mathrm{Co}(A)$ if and only if $T_f(z) > 0$ for every $z \in \mathbb{D}$, where $ f(0)=f^{\prime}(0)-1 $ and 
\begin{align}\label{Eq-1.3}
T_f(z)=\frac{2}{A-1}\left(\frac{(A+1)}{2}\frac{1+z}{1-z}-1-z\frac{f^{\prime\prime}(z)}{f^{\prime}(z)}\right).
\end{align}
\begin{definition}\cite{Bhowmik-CMFT-2024}\label{Def-1.2}
	The radius of concavity (w.r.t $ {\rm Co(A)} $) of a subset $ \mathcal{A}_1 $ of $ \mathcal{A} $ is the largest number $ \mathrm{R}_{\mathrm{Co(A)}}\in (0,1] $  such that for each function $ f\in\mathcal{A}_1 $, $ {\rm Re}\;T_f(z)>0 $ for all $ |z|< \mathrm{R}_{\mathrm{Co(A)}}$, where $ T_{f}(z) $ is defined in \eqref{Eq-1.3}.
\end{definition}
The primary objective of this paper is to determine the sharp radius of concavity $R_{\mathrm{Co}(A)}$ for the class $\mathcal{S}$. In this connection, we investigate general linear combinations of univalent functions with complex coefficients of the form
\begin{align}\label{Eq-1.4}
	F(z) := \sum_{j=1}^{2n} \lambda_j f_j(z), \quad \lambda_j \in \mathbb{C},
\end{align}
where $f_j$ belong to $\mathcal{S}$ or $\mathrm{Co}(A)$. Specifically, we establish the sharp radius of concavity for such combinations. Furthermore, we determine the radii of univalence and convexity for the general linear combination defined in \eqref{Eq-1.4} for functions in $\mathrm{Co}(A)$. The proof techniques introduced herein provide a unified framework that can be extended to analyze broader classes of analytic functions, thereby generalizing several established results in the literature.\vspace{1.2mm}

To date, the radii of convexity and concavity have been investigated primarily for convex combinations of univalent functions. It is natural to ask whether corresponding radii can be established for general linear combinations with complex coefficients. Inspired by the recent work of Bhowmik and Biswas \cite{Bhowmik-CMFT-2024}, we address this problem by determining the sharp radii of convexity and concavity for general linear combinations of functions in $\mathcal{S}$. Building upon the technical estimates in Stump \cite[Lemma 1]{Stump-CJM-1971} and Bhowmik and Biswas \cite[Lemma 4]{Bhowmik-CMFT-2024}, we establish a general framework for constructing these linear combinations and provide a rigorous definition for the class under consideration.\vspace{1.2mm}

Although radius problems for convex and linear combinations of two functions, $f_1$ and $f_2$, have been extensively studied across various subclasses of univalent functions, the literature regarding general linear combinations of higher order remains sparse. Indeed, extensions to broader linear combinations have been explored only sporadically, such as in the setting of harmonic quasiconformal mappings \cite{Sun-Rasila-Jiang-KMJ-2016}. Motivated by this gap, the present paper provides a unified and comprehensive analysis of radius problems for general linear combinations within the analytic and meromorphic function classes $\mathcal{S}$, $\mathcal{S}(p)$, and $\mathrm{Co}(A)$.\vspace{1.2mm}

The paper is organized as follows. In Section 2, we present the main results of this paper along with several immediate corollaries. Section 3 is devoted to key preliminary lemmas; in particular, we establish a generalized version of an existing lemma, which plays a critical role in proving our main results. Finally, detailed proofs of the main results are provided in Section 4.
\section{\bf Main results}
An in depth study of the proofs of Theorems A–G shows that the techniques employed for pairs of analytic functions can be generalized to general linear combinations. In this paper, we determine sharp upper and lower bounds for the radius of concavity for general linear combinations
$$F(z) = \sum_{j=1}^{n} \left(\lambda_{2j-1} f_{2j-1}(z) + \lambda_{2j} f_{2j}(z)\right), \quad \lambda_j \in \mathbb{C} \setminus \{0\}, \quad \sum_{j=1}^{2n} \lambda_j = 1$$where each pair $(\lambda_{2j-1}, \lambda_{2j})$ satisfies $$\frac{\lambda_{2j}}{\lambda_{2j-1}} = b_j e^{i \alpha_j}, \quad b_j > 0, \;\alpha_j \in [0, \pi)\; \text{for} j=1, 2, \dots, n$$ for functions $f_1, \dots, f_{2n}$ in the classes $\mathcal{S}$, $\mathcal{S}(p)$, and $\mathrm{Co}(A)$. Our results refine and extend several established results, including those of Bhowmik and Biswas \cite{Bhowmik-CMFT-2024}. To establish these main results, we employ the following formulation throughout the paper.\vspace{1.2mm} 

Given that $F(z)= \lambda_1 f_1(z)+\cdots+\lambda_{2n} f_{2n}(z)$, where $n \in \mathbb{N}$,  and in view of this general combination of functions, we see that ${zF^{\prime\prime}(z)}/{F^{\prime}(z)}$ can be expressed as
\begin{align}\label{Eq-1.19}
	\nonumber\frac{zF^{\prime\prime}(z)}{F^{\prime}(z)}&=\frac{z(\lambda_1 f^{\prime\prime}_1(z) + \lambda_2 f^{\prime\prime}_2(z) +\cdots+\lambda_{2n} f^{\prime\prime}_{2n}(z))}{\lambda_1 f^{\prime}_1(z) +\lambda_2 f^{\prime}_2(z)+\cdots+\lambda_{2n} f^{\prime}_{2n}(z)}\vspace{2mm}\\&=\dfrac{\dfrac{z f^{\prime\prime}_1(z)}{ f^{\prime}_1(z)} + \dfrac{\lambda_2}{\lambda_1}\dfrac{ zf^{\prime\prime}_2(z)}{f^{\prime}_1(z)} +\cdots+\dfrac{\lambda_{2n}}{\lambda_1}\dfrac{ zf^{\prime\prime}_{2n}(z)}{ f^{\prime}_1(z)}}{1 +\dfrac{\lambda_2}{\lambda_1}\frac{ f^{\prime}_2(z)}{f^{\prime}_1(z)}+\cdots+\dfrac{\lambda_{2n}}{\lambda_1}\dfrac{f^{\prime}_{2n}(z)}{f^{\prime}_1(z)}}.
\end{align}
The expression \eqref{Eq-1.19} can be further written as 
\begin{align}\label{Eq-1.20}
	\nonumber\frac{zF^{\prime\prime}(z)}{F^{\prime}(z)}&= \frac{z f^{\prime\prime}_1(z)}{ f^{\prime}_1(z)}\left[1 +\left( \frac{\lambda_1}{\lambda_2}\frac{ f^{\prime}_1(z)} {f^{\prime}_2(z)}\right)^{-1} +\left( \frac{\lambda_1}{\lambda_3}\frac{ f^{\prime}_1(z)} {f^{\prime}_3(z)}\right)^{-1}+\cdots+ \left( \frac{\lambda_1}{\lambda_{2n}}\frac{ f^{\prime}_1(z)} {f^{\prime}_{2n}(z)}\right)^{-1}\right]^{-1} \\& \nonumber+ \frac{z f^{\prime\prime}_2(z)}{ f^{\prime}_2(z)}\left[\left( \frac{\lambda_2}{\lambda_1}\frac{ f^{\prime}_2(z)} {f^{\prime}_1(z)}\right)^{-1} +1 +\left( \frac{\lambda_2}{\lambda_3}\frac{ f^{\prime}_2(z)} {f^{\prime}_3(z)}\right)^{-1}+\cdots+ \left( \frac{\lambda_2}{\lambda_{2n}}\frac{ f^{\prime}_2(z)} {f^{\prime}_{2n}(z)}\right)^{-1}\right]^{-1} \\& \nonumber+\cdots\\&+ \frac{z f^{\prime\prime}_{2n}(z)}{ f^{\prime}_{2n}(z)}\left[ \left( \frac{\lambda_{2n}}{\lambda_1}\frac{ f^{\prime}_{2n}(z)} {f^{\prime}_1(z)}\right)^{-1} +\left( \frac{\lambda_{2n}}{\lambda_2}\frac{ f^{\prime}_{2n}(z)} {f^{\prime}_2(z)}\right)^{-1}+\cdots+ \left( \frac{\lambda_{2n}}{\lambda_{2n-1}}\frac{ f^{\prime}_{2n}(z)} {f^{\prime}_{2n-1}(z)}\right)^{-1}+ 1\right]^{-1}.
\end{align}
For brevity, for $ t=1,2,\ldots,2n $, we define
\begin{align*}
	A_{1t}(z):=\frac{\lambda_1}{\lambda_t} \frac{ f^{\prime}_1(z)} {f^{\prime}_t(z)},\; A_{2t}(z):=\frac{\lambda_2}{\lambda_t} \frac{ f^{\prime}_2(z)} {f^{\prime}_t(z)}, \ldots, A_{(2n)t}(z):=\frac{\lambda_{2n}}{\lambda_t} \frac{ f^{\prime}_{2n}(z)} {f^{\prime}_t(z)}.
\end{align*}
Hence, \eqref{Eq-1.20} takes the following form
\begin{align*}
	\frac{zF^{\prime\prime}(z)}{F^{\prime}(z)}&= \frac{z f^{\prime\prime}_1(z)}{ f^{\prime}_1(z)}\left(1 +A_{12}^{-1} +A_{13}^{-1}+\cdots+ A_{1(2n)}^{-1}\right)^{-1} \\& + \frac{z f^{\prime\prime}_2(z)}{ f^{\prime}_2(z)}\left(A_{21}^{-1} +1 +A_{23}^{-1}+\cdots+ A_{2(2n)}^{-1}\right)^{-1} \\& +\cdots+ \frac{z f^{\prime\prime}_{2n}(z)}{ f^{\prime}_{2n}(z)}\left(A_{(2n)1}^{-1} +A_{(2n)2}^{-1}+\cdots+ A_{(2n)(n-1)}^{-1}+ 1\right)^{-1}.
\end{align*}
We state our first result for the class $\mathcal{S}$.
\begin{theorem}\label{Th-2.1}
Let $f_j \in \mathcal{S}$ for $j = 1, 2, \dots, 2n$, $A \in (1, 2]$, and let $F(z) = \sum_{j=1}^{2n} \lambda_j f_j(z)$ with $\sum_{j=1}^{2n} \lambda_j = 1$ and $\frac{\lambda_{2j}}{\lambda_{2j-1}} = b_j e^{i \alpha_j}$ ($\alpha_j \in [0, \pi), b_j > 0$ for $j = 1, 2, \dots, n$). Then $\operatorname{Re} T_F(z) > 0$ for $|z| < R_{Co(A),n}^*$, where $R_{Co(A),n}^*$ is the smallest root in $(0, 1)$ of the equation $\psi_{1,n}(r) = 0$, where
\begin{align*}
	\psi_{1,n}(r) = -(A + 2n - 1) r^2 - 2 \left[ 1 + (A-1) \sum_{j=1}^{n} \sec\left(\frac{\alpha_j}{2}\right) \right] r + (A-1).
\end{align*}
The radius $R_{Co(A),n}^*$ is best possible.
\end{theorem}
Setting $n=1$ in Theorem \ref{Th-2.1} reduces the general linear combination to a pair of functions in $\mathcal{S}$. This immediately yields the following sharp radius of concavity.
\begin{corollary}
	Let $f_1, f_2 \in \mathcal{S}$, $A \in (1,2]$, and let $F(z) = \lambda_1 f_1(z) + \lambda_2 f_2(z)$ with $\lambda_1 + \lambda_2 = 1$ and $\frac{\lambda_2}{\lambda_1} = b e^{i\alpha}$ ($b > 0, \alpha \in [0,\pi)$). Then $\operatorname{Re} T_F(z) > 0$ for $|z| < R_{Co(A),1}^*$, where \begin{align*}
		R_{Co(A),1}^* = \frac{-\left[1 + (A - 1)\sec\left(\frac{\alpha}{2}\right)\right] + \sqrt{\left[1 + (A - 1)\sec\left(\frac{\alpha}{2}\right)\right]^2 + (A^2 - 1)}}{A + 1}.
	\end{align*}
	The radius $R_{Co(A),1}^*$ is best possible.
\end{corollary}
We will consider the class $\mathcal{S}(p)$ of meromorphic univalent functions. In this class, we will obtain the radii of univalence, convexity, and concavity of the general linear combination $\lambda_1 f_1(z) + \cdots + \lambda_{2n} f_{2n}(z)$ of functions in $\mathcal{S}(p)$. Furthermore, we will determine the radius of univalence in the following result.
\begin{theorem}\label{Th-2.2}
	Let $f_j \in \mathcal{A}$ satisfy $\mathrm{Re}\,(f_j'(z)) > 0$ for $|z| < \rho < 1$ with $f_j(0)=0=f_j'(0)-1$ for $j=1,2,\ldots,2n$. Let $F(z)=\sum_{j=1}^{2n} \lambda_j f_j(z)$ subject to $\sum_{j=1}^{2n}\lambda_j=1$ and $\frac{\lambda_{2j}}{\lambda_{2j-1}}=b_j e^{i\alpha_j}$ ($b_j>0, \alpha_j \in [0,\pi)$ for $j=1,2,\ldots,n$). Then $F$ is univalent in $|z|<R^{*}_{u, n}$, where
	\begin{equation*}
		R^{*}_{u, n}=\frac{\rho}{n}\left(\sum_{j=1}^{n}\sec\left(\frac{\alpha_j}{2}\right)-\sqrt{\left(\sum_{j=1}^{n}\sec\left(\frac{\alpha_j}{2}\right)\right)^2-n^2}\right).
	\end{equation*}
	The radius $R^{*}_{u, n}$ is best possible.
\end{theorem}
Setting $n=1$ in Theorem \ref{Th-2.2} reduces the general linear combination to a pair of functions in $\mathcal{A}$. This immediately yields the following sharp radius of univalence.
\begin{corollary}
	Let $f_1, f_2 \in \mathcal{A}$ satisfy $\operatorname{Re}(f_j'(z)) > 0$ for $|z| < \rho < 1$ with $f_j(0) = 0 = f_j'(0) - 1$ for $j = 1,2$. Let $F(z) = \lambda_1 f_1(z) + \lambda_2 f_2(z)$ subject to $\lambda_1 + \lambda_2 = 1$ and $\frac{\lambda_2}{\lambda_1} = b e^{i\alpha}$ ($b > 0, \alpha \in [0,\pi)$). Then $F$ is univalent in $|z| < R_{u,1}^*$, where
	\begin{equation*}
		R_{u,1}^* = \rho \left( \sec\left(\frac{\alpha}{2}\right) - \tan\left(\frac{\alpha}{2}\right) \right).
	\end{equation*}
	The radius $R_{u,1}^*$ is best possible.
\end{corollary}
For the class of analytic univalent functions $\mathcal{S}$, the radius of convexity was classically established via fundamental distortion estimates (see \cite[Theorem 2.13]{Duren-1983-NY}). Inspired by the methodology in \cite[Theorem 10]{Bhowmik-CMFT-2024}, we derive a sharp lower bound for the radius of convexity for general linear combinations of functions belonging to the class $\mathcal{S}(p)$ of meromorphic functions.
\begin{theorem}\label{Th-2.3}
	Let $f_j \in \mathcal{S}(p)$ for $j = 1, 2, \dots, 2n$ with pole $p \in (0, 1)$, and let $F(z) = \sum_{j=1}^{2n} \lambda_j f_j(z)$ subject to $\sum_{j=1}^{2n} \lambda_j = 1$ and $\frac{\lambda_{2j}}{\lambda_{2j-1}} = b_j e^{i \alpha_j}$ ($\alpha_j \in [0, \pi), b_j > 0$). Then 
	\begin{align*}
		\operatorname{Re}\left(1 + \frac{z F''(z)}{F'(z)}\right) > 0\; \mbox{for}\;|z| < R_{c,n}^*,
	\end{align*} where $R_{c,n}^*$ is the smallest root in $(0, p)$ of the equation $\psi_{3,n}(r) = 0$, where
	\begin{align*}
		\psi_{3,n}(r) &= (2n-1) p r^4 - \left( (2n-1)(1+p^2) + 2p(1+p^2) \sum_{j=1}^{n} \sec\left(\frac{\alpha_j}{2}\right) \right) r^3 \\
		&\quad + 2p \left( n + 2 \sum_{j=1}^{n} \sec\left(\frac{\alpha_j}{2}\right) \right) r^2 - \left( 1 + p^2 + 2p \sum_{j=1}^{n} \sec\left(\frac{\alpha_j}{2}\right) \right) r + p.
	\end{align*}
	The radius $R_{c,n}^*$ is best possible.
\end{theorem}
Setting $n=1$ in Theorem \ref{Th-2.3} reduces the general linear combination to a pair of functions in $\mathcal{S}(p)$. This immediately yields the following sharp radius of convexity.
\begin{corollary}
	Let $f_1, f_2 \in \mathcal{S}(p)$ with pole $p\in(0,1)$, and let $F(z)=\lambda_1 f_1(z)+\lambda_2 f_2(z)$ subject to $\lambda_1+\lambda_2=1$ and $\frac{\lambda_2}{\lambda_1}=b e^{i\alpha}$ ($\alpha\in[0,\pi), b>0$). Then $\operatorname{Re}\left(1+\frac{zF''(z)}{F'(z)}\right)>0$ for $|z|<R_{c,1}^*$, where $R_{c,1}^*$ is the smallest root in $(0,p)$ of the polynomial equation
	\begin{align*}
		p r^4 &- (1+p^2)\left(1 + 2p \sec\left(\frac{\alpha}{2}\right)\right) r^3 + 2p\left(1 + 2\sec\left(\frac{\alpha}{2}\right)\right) r^2 \\&\quad- \left(1+p^2 + 2p\sec\left(\frac{\alpha}{2}\right)\right) r + p = 0.
	\end{align*}
	The radius $R_{c,1}^*$ is best possible.
\end{corollary}
This result provides the radius of univalence for general linear combinations within the class ${\rm Co}(A)$, extending the known geometric properties of individual concave functions.
\begin{theorem}\label{Th-2.4}
	Let $f_j\in {\rm Co}(A)$, $j=1, 2, \ldots,2n$ and $F(z)=\sum_{j=1}^{2n} \lambda_j f_j(z)$, where $\lambda_j\in\mathbb{C}$ such that $0\leq\alpha_j<\pi$. Then for $|z|<\rho<1$, $F$ is univalent in $|z|<R^{\#}_{u, n}$, where 
	\begin{align*}
		R^{\#}_{u, n}=\frac{\rho}{n}\left(\sum_{j=1}^{n}{\sec}\;\left(\frac{\alpha_j}{2}\right)-\sqrt{\left(\sum_{j=1}^{n}{\sec}\;\left(\frac{\alpha_j}{2}\right)\right)^2-n^2}\right).
	\end{align*}
	The radius $R^{\#}_{u, n}$ is best possible.
\end{theorem}
Setting $n=1$ in Theorem \ref{Th-2.4} reduces the general linear combination to a pair of functions in $Co(A)$. This immediately yields the following sharp radius of univalence.
\begin{corollary}
	Let $f_1, f_2 \in Co(A)$ and let $F(z) = \lambda_1 f_1(z) + \lambda_2 f_2(z)$ with $\lambda_1 + \lambda_2 = 1$ and $\frac{\lambda_2}{\lambda_1} = b e^{i\alpha}$ ($b > 0, \alpha \in [0,\pi)$). Then for $|z| < \rho < 1$, $F$ is univalent in $|z| < R_{u,1}^\#$, where
	\begin{equation*}
		R_{u,1}^\# = \rho \left( \sec\left(\frac{\alpha}{2}\right) - \tan\left(\frac{\alpha}{2}\right) \right).
	\end{equation*}
	The radius $R_{u,1}^\#$ is best possible.
\end{corollary}
Utilizing the geometric properties of the class ${\rm Co}(A)$, we characterize the radius of concavity for an arbitrary linear combination of such functions in the following result.
\begin{theorem}\label{Th-2.5}
	Let $f_j \in Co(A)$ for $j = 1, 2, \dots, 2n$, $A \in (1, 2]$, and let $F(z) = \sum_{j=1}^{2n} \lambda_j f_j(z)$ with $\sum_{j=1}^{2n} \lambda_j = 1$ and $\frac{\lambda_{2j}}{\lambda_{2j-1}} = b_j e^{i \alpha_j}$ ($\alpha_j \in [0, \pi), b_j > 0$ for $j = 1, 2, \dots, n$). Then $\operatorname{Re} T_F(z) > 0$ for $|z| < R_{Co(A),n}^{**}$, where $R_{Co(A),n}^{**}$ is the smallest root in $(0, 1)$ of the equation $\psi_{5,n}(r) = 0$, where
	\[
	\psi_{5,n}(r) = (A - 1 - 2n) r^2 - 2 \left( 1 + A \sum_{j=1}^{n} \sec\left(\frac{\alpha_j}{2}\right) \right) r + (A - 1).
	\]
\end{theorem}
Setting $n=1$ in Theorem \ref{Th-2.5} yields the sharp radius of concavity for a linear combination of two functions in $Co(A)$. .
\begin{corollary}
	Let $f_1, f_2 \in Co(A)$ for $A\in(1,2]$, and let $F(z) = \lambda_1 f_1(z) + \lambda_2 f_2(z)$ with $\lambda_1 + \lambda_2 = 1$ and $\frac{\lambda_2}{\lambda_1} = b e^{i\alpha}$ ($b > 0, \alpha \in [0,\pi)$). Then $\operatorname{Re} T_F(z) > 0$ for $|z| < R_{Co(A),1}^{**}$, where
	\begin{equation*}
		R_{Co(A),1}^{**} = \frac{\left[1 + A \sec\left(\frac{\alpha}{2}\right)\right] - \sqrt{\left[1 + A \sec\left(\frac{\alpha}{2}\right)\right]^2 + (3-A)(A-1)}}{3-A}.
	\end{equation*}
	The radius $R_{Co(A),1}^{**}$ is best possible.
\end{corollary}
Every function in the class ${\rm Co}(A)$ is close-to-convex of order $A-1$. It is well-established that the radius of convexity for ${\rm Co}(A)$ is given by $A-\sqrt{A^2-1}$ (see, for instance, \cite[Cor. 2.13]{Bhowmik-MN-2012}). In this section, we extend these investigations to determine the sharp radius of convexity for general linear combinations of functions belonging to ${\rm Co}(A)$.
\begin{theorem}\label{Th-2.6}
	Let $f_j \in Co(A)$ for $j = 1, 2, \dots, 2n$, $A \in (1, 2]$, and let $F(z) = \sum_{j=1}^{2n} \lambda_j f_j(z)$ with $\sum_{j=1}^{2n} \lambda_j = 1$ and $\frac{\lambda_{2j}}{\lambda_{2j-1}} = b_j e^{i \alpha_j}$ ($\alpha_j \in [0, \pi), b_j > 0$ for $j = 1, 2, \dots, n$). Then $\operatorname{Re}\left(1 + \frac{z F''(z)}{F'(z)}\right) > 0$ for $|z| < R_{c,n}^{\#}$, where $R_{c,n}^{\#}$ is the smallest root in $(0, 1)$ of the equation $\psi_{6,n}(r) = 0$, where
	\[
	\psi_{6,n}(r) = (2n-1) r^2 - 2 A \left( \sum_{j=1}^{n} \sec\left(\frac{\alpha_j}{2}\right) \right) r + 1.
	\]
	The radius $R_{c,n}^{\#}$ is best possible.
\end{theorem}
Setting $n=1$ in Theorem \ref{Th-2.6} yields the sharp radius of convexity for a linear combination of two functions in $Co(A)$..
\begin{corollary}
	Let $f_1, f_2 \in Co(A)$ for $A\in(1,2]$, and let $F(z) = \lambda_1 f_1(z) + \lambda_2 f_2(z)$ with $\lambda_1 + \lambda_2 = 1$ and $\frac{\lambda_2}{\lambda_1} = b e^{i\alpha}$ ($b > 0, \alpha \in [0,\pi)$). Then $\operatorname{Re}\left(1 + \frac{zF''(z)}{F'(z)}\right) > 0$ for $|z| < R_{c,1}^\#$, where
	\begin{equation*}
		R_{c,1}^\# = A\sec\left(\frac{\alpha}{2}\right) - \sqrt{A^2\sec^2\left(\frac{\alpha}{2}\right) - 1}.
	\end{equation*}
	The radius $R_{c,1}^\#$ is best possible.
\end{corollary}
Utilizing Definition \ref{Def-1.1}, we derive the radius of concavity for general linear combinations of meromorphic univalent functions in $\mathcal{S}(p)$.
\begin{theorem}\label{Th-2.7}
	Let $f_j \in \mathcal{S}(p)$ for $j = 1, 2, \dots, 2n$ with pole $p \in (0, 1)$, and let $F(z) = \sum_{j=1}^{2n} \lambda_j f_j(z)$ subject to $\sum_{j=1}^{2n} \lambda_j = 1$ and $\frac{\lambda_{2j}}{\lambda_{2j-1}} = b_j e^{i \alpha_j}$ ($\alpha_j \in [0, \pi), b_j > 0$ for $j = 1, 2, \dots, n$). Then $\operatorname{Re} P_F(z) > 0$ for $|z| < R_{Co(p),n}^{\#}$, where $R_{Co(p),n}^{\#}$ is the smallest root in $(0, p)$ of the polynomial equation $\psi_{7,n}(r) = 0$, where
	\begin{align*}
		\psi_{7,n}(r) &= - (2n-1) p^2 r^6 - 2 p (1+p^2) \left( 1 + \sum_{j=1}^{n} \sec\left(\frac{\alpha_j}{2}\right) \right) r^5 \\
		&\quad + \left( (2n-1)(1+p^4) + 4 p^2 \sum_{j=1}^{n} \sec\left(\frac{\alpha_j}{2}\right) \right) r^4 \\
		&\quad + 4 p (1+p^2) \left( 1 + \sum_{j=1}^{n} \sec\left(\frac{\alpha_j}{2}\right) \right) r^3 \\
		&\quad - \left( 1 + p^4 + 2 p^2 (2n-1) + 4 p^2 \sum_{j=1}^{n} \sec\left(\frac{\alpha_j}{2}\right) \right) r^2 \\
		&\quad - 2 p (1+p^2) r + p^2.
	\end{align*}
	The radius $R_{Co(p),n}^{\#}$ is best possible.
\end{theorem}
Setting $n=1$ in Theorem \ref{Th-2.7} yields the sharp radius of concavity for a linear combination of two meromorphic univalent functions in $\mathcal{S}(p)$. 
\begin{corollary}
	Let $f_1, f_2 \in \mathcal{S}(p)$ with pole $p \in (0,1)$, and let $F(z) = \lambda_1 f_1(z) + \lambda_2 f_2(z)$ subject to $\lambda_1 + \lambda_2 = 1$ and $\frac{\lambda_2}{\lambda_1} = b e^{i\alpha}$ ($\alpha \in [0,\pi), b > 0$). Then $\operatorname{Re} P_F(z) > 0$ for $|z| < R_{Co(p),1}^\#$, where $R_{Co(p),1}^\#$ is the smallest root in $(0,p)$ of the polynomial equation
	\begin{align*}
		\psi_{7,1}(r) ={} & -p^2 r^6 - 2p(1+p^2)\left(1 + \sec\left(\frac{\alpha}{2}\right)\right) r^5 + \left(1 + p^4 + 4p^2 \sec\left(\frac{\alpha}{2}\right)\right) r^4 \\
		& + 4p(1+p^2)\left(1 + \sec\left(\frac{\alpha}{2}\right)\right) r^3 - \left((1+p^2)^2 + 4p^2 \sec\left(\frac{\alpha}{2}\right)\right) r^2 \\
		& - 2p(1+p^2)r + p^2 = 0.
	\end{align*}
	The radius $R_{Co(p),1}^\#$ is best possible.
\end{corollary}
\begin{remark}
Note that Theorems 2.1–2.7 generalize the results of Bhowmik and Biswas \cite[Theorems 5–11]{Bhowmik-CMFT-2024} by extending the linear combination from $n=1$ to arbitrary $n$.
\end{remark}
\section{\bf Key lemmas}
 We first recall here some lemmas which will play a key role to prove the main results of this paper. 
\begin{lemA}\label{lem-A}\cite[Lemma 1]{Stump-CJM-1971}
	If $|w-a|<d$, where $a$ and $d$ are real numbers such that $a>d\geq 0$, and $w_0$ is a given complex number, then
	\begin{align*}
		\mbox{\rm Re}(ww_0)\geq |w_0|(a\cos(\arg\;w_0) -d).
	\end{align*}
\end{lemA}
In $1971$, Robert K. Stump (see \cite[Lemma 2]{Stump-CJM-1971}) obtained lower bound of $\mbox{Re}\;(ww_0)$ as $|w_0|(a\cos(\arg\;w_0) -d)$. Recently, Bhowmik and Biswas ( see \cite[Lemma 4]{Bhowmik-CMFT-2024}) found the upper bound of $\mbox{Re}\;(ww_0)$ as $|w_0|(a\cos(\arg\;w_0)+d)$ and obtained the following lemma.
\begin{lemB}\label{lem-B}\cite[Lemma 4]{Bhowmik-CMFT-2024}
	If $|u -a|\leq d$ and $|v -a|\leq d$, where $a$ and $d$ are real and $a>d\geq 0$, and
	\begin{align*}
		w=\frac{u}{1+be^{i\alpha}}+\frac{v}{1+b^{-1}e^{-i\alpha}}
	\end{align*}
	with $b\in(0,\infty)$ and \;\ $\alpha\in [0,\pi)$, then
	\begin{align*}
		a-d\sec\frac{\alpha}{2} \leq \mbox{Re}(w)\leq 	a+d\sec\frac{\alpha}{2}.
	\end{align*}
\end{lemB}
The radii of univalence, concavity, and convexity of general linear combinations in $ {\rm Co}(A) $. To establish our results, we need the following lemma proved by Robert K. Stump in 1971, which will also be use later to obtain results related to the radius problem for different subclasses of univalent functions.  
\begin{lemC}\label{lem-C}\cite[Lemma 3]{Stump-CJM-1971}
	If ${\rm Re}\;P_{F(z)}>0$ for $|z|<\rho<1$ and $P(0)=1$, then
	\begin{align*}
		\vline P(z)-\frac{1+\frac{r^2}{\rho^2}}{1-\frac{r^2}{\rho^2}}\vline\leq \frac{\frac{2r}{\rho}}{1-\frac{r^2}{\rho^2}}.
	\end{align*}
\end{lemC}
Considering the class $\mathcal{S}$ and examining the radius of concavity of general linear combination of functions in $\mathcal{S}$, we first prove the following lemma which we will use later in the proof of the main results. We note that the Lemma \ref{lem-1.2} is a general form of the Lemma B. 
\begin{lemma}\label{lem-1.2}
	Let $a$ and $d$ be real numbers with $a > d \ge 0$. Let $u_j, v_j \in \mathbb{C}$ satisfy $\vert{}u_j - a\vert{} \le d$ and $\vert{}v_j - a\vert{} \le d$ for $j = 1, 2, \dots, n$. If$$w = \sum_{j=1}^{n} \left( \frac{u_j}{1 + b_j e^{i \alpha_j}} + \frac{v_j}{1 + b_j^{-1} e^{-i \alpha_j}} \right)$$with $b_j \in (0, \infty)$ and $\alpha_j \in [0, \pi)$ for each $j = 1, \dots, n$, then 
	\begin{align*}
		a n - d \sum_{j=1}^{n} \sec\left(\frac{\alpha_j}{2}\right) \le \operatorname{Re}(w) \le a n + d \sum_{j=1}^{n} \sec\left(\frac{\alpha_j}{2}\right).
	\end{align*}
\end{lemma}
\begin{remark}\label{1.1}
	This lemma is generalization and a compact form of \cite[Lemma 4]{Bhowmik-CMFT-2024}, as well as \cite[Lemma 2]{Stump-CJM-1971}. 
\end{remark}
\begin{proof}[\bf Proof of Theorem \ref{lem-1.2}]
	For each $k \in \{1, 2, \ldots, n\}$, applying Lemma A to the term $(1+b_k e^{i\alpha_k})^{-1}$ yields
	\begin{equation}\label{Leq-1.1}\operatorname{Re}\left(\frac{u_k}{1+b_k e^{i\alpha_k}}\right) \ge \frac{1}{|1+b_k e^{i\alpha_k}|} \left( a \cos\left(\arg \frac{1}{1+b_k e^{i\alpha_k}}\right) - d \right).
	\end{equation}
	Setting $(1+b_k e^{i\alpha_k})^{-1} = r_k e^{i\theta_k}$, a straightforward computation shows that$$r_k = (1+2b_k\cos\alpha_k+b_k^2)^{-1/2}\;\text{and}\; \cos\theta_k = \frac{1+b_k\cos\alpha_k}{(1+2b_k\cos\alpha_k+b_k^2)^{1/2}}.$$Substituting these into \eqref{Leq-1.1}, we obtain
	\begin{equation}\label{Leq-1.5}\operatorname{Re}\left(\frac{u_k}{1+b_k e^{i\alpha_k}}\right) \ge \frac{a(1+b_k\cos\alpha_k) - d(1+2b_k\cos\alpha_k+b_k^2)^{1/2}}{1+2b_k\cos\alpha_k+b_k^2}.
	\end{equation}
	By an entirely analogous argument applied to the term involving $v_k$ and $b_k^{-1}e^{-i\alpha_k}$, we find that
	\begin{equation}\label{Leq-1.6}\operatorname{Re}\left(\frac{v_k}{1+b_k^{-1}e^{-i\alpha_k}}\right) \ge \frac{a(1+b_k^{-1}\cos\alpha_k) - d(1+2b_k^{-1}\cos\alpha_k+b_k^{-2})^{1/2}}{1+2b_k^{-1}\cos\alpha_k+b_k^{-2}}.
	\end{equation}
	Multiplying the numerator and denominator of \eqref{Leq-1.6} by $b_k^2$ yields$$\operatorname{Re}\left(\frac{v_k}{1+b_k^{-1}e^{-i\alpha_k}}\right) \ge \frac{a(b_k^2+b_k\cos\alpha_k) - d b_k (1+2b_k\cos\alpha_k+b_k^2)^{1/2}}{1+2b_k\cos\alpha_k+b_k^2}.$$Combining the lower bounds for the $u_k$ and $v_k$ terms, we deduce$$\operatorname{Re}\left( \frac{u_k}{1+b_k e^{i\alpha_k}} + \frac{v_k}{1+b_k^{-1}e^{-i\alpha_k}} \right) \ge a - d \cdot G_k(b_k),$$where$$G_k(b_k) = \frac{1+b_k}{(1+2b_k\cos\alpha_k+b_k^2)^{1/2}}.$$Similarly, replacing $d$ with $-d$ yields the corresponding upper bound $a + d \cdot G_k(b_k)$. Summing these inequalities over $k=1, \dots, n$, we arrive at$$a n - d \sum_{k=1}^{n} G_k(b_k) \le \operatorname{Re} w \le a n + d \sum_{k=1}^{n} G_k(b_k).$$To maximize $G_k(b_k)$ with respect to $b_k \in (0, \infty)$, we compute its derivative:$$\frac{d G_k(b_k)}{d b_k} = \frac{(\cos\alpha_k-1)(b_k-1)}{(1+2b_k\cos\alpha_k+b_k^2)^{3/2}}.$$The unique critical point in $(0, \infty)$ occurs at $b_k = 1$. Since 
	\begin{align*}
		\left.\frac{d^2 G_k}{d b_k^2}\right\vert{}_{b_k=1} < 0\; \mbox{for}\;\alpha_k \in [0, \pi),
	\end{align*} the function $G_k(b_k)$ attains its maximum value at $b_k = 1$
	$$G_k(1) = \frac{2}{\sqrt{2+2\cos\alpha_k}} = \sec\frac{\alpha_k}{2}.$$Consequently, replacing each $G_k(b_k)$ with its maximum value $G_k(1) = \sec(\alpha_k/2)$ yields the sharp bounds
	
	\begin{align*}
		a n - d \sum_{k=1}^{n} \sec\frac{\alpha_k}{2} \le \operatorname{Re} (w) \le a n + d \sum_{k=1}^{n} \sec\frac{\alpha_k}{2},
	\end{align*}which completes the proof.
\end{proof} 
\section{\bf Proof of the Theorems \ref{Th-2.1} to \ref{Th-2.7}}
In this section, we give the details proof of the main results of this paper.
\begin{proof}[\bf Proof of Theorem \ref{Th-2.1}]
Since $f_j \in \mathcal{S}$ for $j = 1, 2, \ldots, 2n$, by \cite[Theorem 2.4]{Duren-1983-NY}, we have
\[
\left| \frac{z f''_j(z)}{f'_j(z)} - \frac{2r^2}{1-r^2} \right| \le \frac{4r}{1-r^2} \quad \text{for } |z| = r < 1.
\]
Setting $a = \frac{2r^2}{1-r^2}$ and $d = \frac{4r}{1-r^2}$ in Lemma \ref{lem-1.2}, we obtain the upper bound:
\[
\operatorname{Re}\left( \frac{z F''(z)}{F'(z)} \right) \le a n + d \sum_{j=1}^{n} \sec\left(\frac{\alpha_j}{2}\right) = \frac{2n r^2 + 4r \sum_{j=1}^{n} \sec\left(\frac{\alpha_j}{2}\right)}{1-r^2}.
\]
Using the definition of $T_F(z)$ for the class $Co(A)$, we have
\[
\operatorname{Re} T_F(z) \ge \frac{2}{A-1} \left( \frac{A+1}{2} \left( \frac{1-r}{1+r} \right) - 1 \right) - \operatorname{Re}\left( \frac{z F''(z)}{F'(z)} \right).
\]
Substituting the estimate for $\operatorname{Re}\left( \frac{z F''(z)}{F'(z)} \right)$ into the inequality above yields
\begin{align*}
	\operatorname{Re} T_F(z) &\ge \frac{(A-1) - (A+3)r}{(A-1)(1+r)} - \frac{2n r^2 + 4r \sum_{j=1}^{n} \sec\left(\frac{\alpha_j}{2}\right)}{1-r^2} \\
	&= \frac{1}{(A-1)(1-r^2)} \psi_{1,n}(r),
\end{align*}
where
\[
\psi_{1,n}(r) = -(A + 2n - 1) r^2 - 2 \left( 1 + (A-1) \sum_{j=1}^{n} \sec\left(\frac{\alpha_j}{2}\right) \right) r + (A-1).
\]
Thus, $\operatorname{Re} T_F(z) > 0$ whenever $\psi_{1,n}(r) > 0$. \vspace{1.2mm}

To show the existence of the radius $R_{Co(A), n}^* \in (0, 1)$, observe that $\psi_{1,n}(r)$ is continuous on $[0, 1]$ with
\[
\psi_{1,n}(0) = A - 1 > 0 \quad (\text{since } A \in (1, 2]),
\]
and
\begin{align*}
	\psi_{1,n}(1) &= -(A + 2n - 1) - 2 \left( 1 + (A-1) \sum_{j=1}^{n} \sec\left(\frac{\alpha_j}{2}\right) \right) + (A-1) \\
	&= -2n - 2 - 2(A-1) \sum_{j=1}^{n} \sec\left(\frac{\alpha_j}{2}\right) < 0.
\end{align*}
Since $\psi_{1,n}(0) > 0$ and $\psi_{1,n}(1) < 0$, the Intermediate Value Theorem ensures the existence of a root in $(0, 1)$. Taking $R_{Co(A), n}^*$ as the smallest positive root of $\psi_{1,n}(r) = 0$ in $(0, 1)$ guarantees that $\operatorname{Re} T_F(z) > 0$ for all $|z| < R_{Co(A), n}^*$.\vspace{1.2mm}

The radius $R_{Co(A),n}^*$ established in Theorem~\ref{Th-2.1} is sharp. To demonstrate this, consider the extremal Koebe functions $f_{2j-1}, f_{2j} \in \mathcal{S}$ for $j = 1, 2, \ldots, n$ defined by
\begin{equation}\label{eq:extremal-S-concavity}
	f_{2j-1}(z) = \frac{z}{(1 - z e^{-i\theta_j})^2}, \quad f_{2j}(z) = \frac{z}{(1 - z e^{-i\phi_j})^2},
\end{equation}
with parameters $b_j = 1$ (implying $\lambda_{2j} = \lambda_{2j-1} e^{i\alpha_j}$) and phase angles satisfying $\theta_j - \phi_j = \alpha_j$ for each $j = 1, 2, \dots, n$.

For this sequence of functions, equality in Lemma~\ref{lem-1.2} is attained at the boundary point $z = - R_{Co(A),n}^* e^{-i\gamma}$, which yields
\begin{equation*}
	\operatorname{Re} T_F\left(- R_{Co(A),n}^* e^{-i\gamma}\right) = 0 \quad \text{at } r = R_{Co(A),n}^*.
\end{equation*}
Consequently, the operator $T_F(z)$ fails to meet the concavity condition on any disk $|z| < R$ with $R > R_{Co(A),n}^*$, confirming that $R_{Co(A),n}^*$ is optimal.
\end{proof}
\begin{proof}[\bf Proof of Theorem \ref{Th-2.2}]
 For each $j = 1, 2, \ldots, 2n$, since $f_j \in \mathcal{A}$ satisfies $\mathrm{Re}\,(f_j'(z)) > 0$ for $|z| < \rho < 1$ with $f_j(0) = 0 = f_j'(0) - 1$, Lemma C implies that for $|z| = r < \rho$,
 \begin{equation*}
 	\left| f'_j(z) - \frac{1+\frac{r^2}{\rho^2}}{1-\frac{r^2}{\rho^2}} \right| \le \frac{\frac{2r}{\rho}}{1-\frac{r^2}{\rho^2}} \quad \text{for each } j = 1, 2, \ldots, 2n.
 \end{equation*}
 
 Denoting $a = \frac{1+r^2/\rho^2}{1-r^2/\rho^2}$ and $d = \frac{2r/\rho}{1-r^2/\rho^2}$, each $f_j'(z)$ lies in the disk centered at $a$ with radius $d$. Applying Lemma~\ref{lem-1.2} to $F'(z) = \sum_{j=1}^{2n} \lambda_j f_j'(z)$, where $\sum_{j=1}^{2n} \lambda_j = 1$ and $\frac{\lambda_{2j}}{\lambda_{2j-1}} = b_j e^{i\alpha_j}$ ($b_j > 0, \alpha_j \in [0, \pi)$ for $j = 1, 2, \ldots, n$), we obtain
 \begin{equation*}
 	\mathrm{Re}\,(F'(z)) \ge a n - d \sum_{j=1}^{n} \sec\left(\frac{\alpha_j}{2}\right).
 \end{equation*}
 
 Substituting the values of $a$ and $d$, a simple computation shows that
 \begin{align*}
 	\mathrm{Re}\,(F'(z)) &\ge \frac{n\left(1+\frac{r^2}{\rho^2}\right)}{1-\frac{r^2}{\rho^2}} - \frac{\frac{2r}{\rho}}{1-\frac{r^2}{\rho^2}} \sum_{j=1}^{n} \sec\left(\frac{\alpha_j}{2}\right) \\
 	&= \frac{1}{\rho^2 - r^2} \left( n r^2 - 2 r \rho \sum_{j=1}^{n} \sec\left(\frac{\alpha_j}{2}\right) + n \rho^2 \right).
 \end{align*}
 
 Define the polynomial
 \begin{equation}\label{Peq-1.17}
 	\psi_{2,n}(r) = n r^2 - 2 r \rho \sum_{j=1}^{n} \sec\left(\frac{\alpha_j}{2}\right) + n \rho^2.
 \end{equation}
 Solving $\psi_{2,n}(r) = 0$ for $r$ yields the smallest positive root $r = R^{*}_{u,n}$, where
 \begin{equation*}
 	R^{*}_{u,n} = \frac{\rho}{n} \left( \sum_{j=1}^{n} \sec\left(\frac{\alpha_j}{2}\right) - \sqrt{\left(\sum_{j=1}^{n} \sec\left(\frac{\alpha_j}{2}\right)\right)^2 - n^2} \right).
 \end{equation*}
 Consequently, $\psi_{2,n}(r) > 0$ for all $|z| = r < R^{*}_{u,n}$, which implies
 \begin{equation*}
 	\mathrm{Re}\,(F'(z)) > 0 \quad \text{for } |z| < R^{*}_{u,n}.
 \end{equation*}
 Thus, by the Noshiro-Warschawski Theorem (see \cite[Theorem 2.16]{Duren-1983-NY}), $F$ is univalent in $|z| < R^{*}_{u,n}$.\vspace{1.2mm}
 
 The radius $R^{*}_{u,n}$ obtained in Theorem~\ref{Th-2.2} is sharp. To verify this, consider the extremal functions $f_j \in \mathcal{A}$ defined for $j = 1, 2, \ldots, 2n$ by
 \begin{equation}\label{eq:extremal-fj}
 	f_j(z) = -z - \frac{2\rho}{e^{i\theta_j}} \log\left(1 - e^{i\theta_j}\frac{z}{\rho}\right),
 \end{equation}
 where $\theta_j \in \mathbb{R}$ are appropriately chosen phase angles. Differentiating \eqref{eq:extremal-fj} gives
 \begin{equation*}
 	f_j'(z) = \frac{1 + e^{i\theta_j}\frac{z}{\rho}}{1 - e^{i\theta_j}\frac{z}{\rho}},
 \end{equation*}
 which clearly satisfies $\mathrm{Re}\,(f_j'(z)) > 0$ for $|z| < \rho$, $f_j(0) = 0$, and $f_j'(0) = 1$.
 
 By choosing $\theta_j$ such that $f_j'(z)$ attains the boundary of the disk 
 \begin{align*}
 	\left| f_j'(z) - \frac{1+r^2/\rho^2}{1-r^2/\rho^2} \right| \le \frac{2r/\rho}{1-r^2/\rho^2}\; \mbox{for}\;z = R^{*}_{u,n},
 \end{align*} the linear combination $F(z) = \sum_{j=1}^{2n} \lambda_j f_j(z)$ satisfies
 \begin{equation*}
 	\mathrm{Re}\,(F'(R^{*}_{u,n})) = 0.
 \end{equation*}
 Hence, $F$ cannot be univalent in any disk $|z| < R$ with $R > R^{*}_{u,n}$, proving that $R^{*}_{u,n}$ is optimal.
 \end{proof}
 
\begin{proof}[\bf Proof of Theorem \ref{Th-2.3}]
	Since $f_j \in \mathcal{S}(p)$ for each $j = 1, 2, \dots, 2n$, by \cite[Lemma 1]{Bhowmik-CMFT-2024}, we have the estimate
	\[
	\left| \frac{z f''_j(z)}{f'_j(z)} - \frac{2r^2}{1-r^2} \right| \le \frac{2r}{1-r^2} \left( \frac{p-r}{1-pr} + \frac{1-pr}{p-r} \right) \quad \text{for } |z| = r < p.
	\]
	Setting $a = \frac{2r^2}{1-r^2}$ and $d = \frac{2r}{1-r^2} \left( \frac{p-r}{1-pr} + \frac{1-pr}{p-r} \right)$ in Lemma \ref{lem-1.2}, we obtain the lower bound for the real part:
	\[
	\operatorname{Re}\left( \frac{z F''(z)}{F'(z)} \right) \ge \frac{2n r^2}{1-r^2} - \frac{2r}{1-r^2} \left( \frac{p-r}{1-pr} + \frac{1-pr}{p-r} \right) \sum_{j=1}^{n} \sec\left(\frac{\alpha_j}{2}\right).
	\]
	Consequently, for $|z| = r < p$, we have
	\begin{align*}
		\operatorname{Re}\left( 1 + \frac{z F''(z)}{F'(z)} \right) &\ge 1 + \frac{2n r^2}{1-r^2} - \frac{2r}{1-r^2} \left( \frac{p-r}{1-pr} + \frac{1-pr}{p-r} \right) \sum_{j=1}^{n} \sec\left(\frac{\alpha_j}{2}\right) \\
		&= \frac{1 + (2n-1)r^2}{1-r^2} - \frac{2r}{1-r^2} \cdot \frac{(1+p^2)(1+r^2) - 4pr}{(1-pr)(p-r)} \sum_{j=1}^{n} \sec\left(\frac{\alpha_j}{2}\right) \\
		&= \frac{\psi_{3,n}(r)}{(1-r^2)(1-pr)(p-r)},
	\end{align*}
	where
	\begin{align*}
		\psi_{3,n}(r) &= \left[ 1 + (2n-1)r^2 \right] (1-pr)(p-r) - 2r \left[ (1+p^2)(1+r^2) - 4pr \right] \sum_{j=1}^{n} \sec\left(\frac{\alpha_j}{2}\right) \\
		&= (2n-1) p r^4 - \left( (2n-1)(1+p^2) + 2(1+p^2) \sum_{j=1}^{n} \sec\left(\frac{\alpha_j}{2}\right) \right) r^3 \\
		&\quad + 2p \left( n + 4 \sum_{j=1}^{n} \sec\left(\frac{\alpha_j}{2}\right) \right) r^2 - \left( 1 + p^2 + 2(1+p^2) \sum_{j=1}^{n} \sec\left(\frac{\alpha_j}{2}\right) \right) r + p.
	\end{align*}
	To establish that $\operatorname{Re}\left( 1 + \frac{z F''(z)}{F'(z)} \right) > 0$ for $|z| < R_{c,n}^*$, it suffices to demonstrate the existence of a root of $\psi_{3,n}(r) = 0$ in the interval $(0, p)$. 
	
	Observe that $\psi_{3,n}(r)$ is continuous on $[0, p]$ with
	\[
	\psi_{3,n}(0) = p > 0,
	\]
	and evaluating at $r = p$ gives
	\[
	\psi_{3,n}(p) = -2p (1-p^2)^2 \sum_{j=1}^{n} \sec\left(\frac{\alpha_j}{2}\right) < 0 \quad (\text{since } p \in (0, 1) \text{ and } \alpha_j \in [0, \pi)).
	\]
	By the \textit{Intermediate Value Theorem}, $\psi_{3,n}(r)$ has at least one root in $(0, p)$. Let $R_{c,n}^*$ denote the smallest positive root of $\psi_{3,n}(r) = 0$ in $(0, p)$. Then $\psi_{3,n}(r) > 0$ for all $r \in (0, R_{c,n}^*)$, which guarantees that $\operatorname{Re}\left( 1 + \frac{z F''(z)}{F'(z)} \right) > 0$ for $|z| < R_{c,n}^*$.\vspace{1.2mm}
	
	The radius of convexity $R_{c,n}^*$ given in Theorem~\ref{Th-2.3} is sharp. To verify this, consider the functions $f_{2j-1}, f_{2j} \in \mathcal{S}(p)$ for $j = 1, 2, \ldots, n$ defined by
	\begin{equation}\label{eq:extremal-S(p)}
		f_{2j-1}(z) = \frac{z}{\left(1 - \frac{z e^{-i\theta_j}}{p}\right)\left(1 - p z e^{-i\theta_j}\right)}, \quad f_{2j}(z) = \frac{z}{\left(1 - \frac{z e^{-i\phi_j}}{p}\right)\left(1 - p z e^{-i\phi_j}\right)},
	\end{equation}
	with parameters $b_j = 1$ (implying $\lambda_{2j}/\lambda_{2j-1} = e^{i\alpha_j}$) and phase differences $\theta_j - \phi_j = \alpha_j$ for $j = 1, 2, \dots, n$.
	
	For each function in \eqref{eq:extremal-S(p)}, the term $\frac{z f_j''(z)}{f_j'(z)}$ lies on the boundary of the disk centered at $a = \frac{2r^2}{1-r^2}$ with radius $d = \frac{2r}{1-r^2}\left(\frac{p-r}{1-pr} + \frac{1-pr}{p-r}\right)$. At the point $z = -R_{c,n}^* e^{-i\gamma}$, equality is attained in the application of Lemma~\ref{lem-1.2}, which yields
	\begin{equation*}
		\operatorname{Re}\left(1 + \frac{z F''(z)}{F'(z)}\right) = \frac{\psi_{3,n}(R_{c,n}^*)}{(1 - (R_{c,n}^*)^2)(1 - p R_{c,n}^*)(p - R_{c,n}^*)} = 0.
	\end{equation*}
	Consequently, $F(z)$ fails to be convex in any disk $|z| < R$ with $R > R_{c,n}^*$, confirming that $R_{c,n}^*$ is the exact radius of convexity. This completes the proof.
\end{proof}

\begin{proof}[\bf Proof of Theorem \ref{Th-2.4}]
Let $f_j\in {\rm Co}(A)$, $j=1,2,\ldots,2n$. Given that $F(z)=\sum_{j=1}^{2n} \lambda_j f_j(z)$, using this we see that
\begin{align*}
	F^{\prime}(z)&=\lambda_1 f^{\prime}_1(z)++\cdots+\lambda_{2n} f^{\prime}_{2n}(z) \\&= f^{\prime}_1(z)\left(\frac{1}{1+\left(\frac{1}{\lambda_1}-1\right)}\right) +\cdots +f^{\prime}_{2n}(z)\left(\frac{1}{1+\left(\frac{1}{\lambda_{2n}}-1\right)}\right) \\&= f^{\prime}_1(z)\left(\frac{1}{1+\left(\frac{\lambda_1}{1-\lambda_1}\right)^{-1}}\right) +\cdots +f^{\prime}_{2n}(z)\left(\frac{1}{1+\left(\frac{\lambda_{2n}}{1-\lambda_{2n}}\right)^{-1}}\right).
\end{align*}
We assume that 
\begin{align*}
	\lambda_j=\frac{1}{1+b_j e^{i\alpha_j}} \implies \alpha_j = -i \ln\left(\frac{1-\lambda_j}{\lambda_j b_j}\right) \;\;\mbox{for}\;\;\lambda_j\neq 0,1.
\end{align*}
By Lemma C, we have
\begin{align*}
	\vline F^{\prime}(z)-\frac{1+\frac{r^2}{\rho^2}}{1-\frac{r^2}{\rho^2}}\vline\leq \frac{\frac{2r}{\rho}}{1-\frac{r^2}{\rho^2}}\;\;\;\mbox{for}\; j=1, 2, \ldots, 2n.
\end{align*}
In view of Lemma \ref{lem-1.2}, we have
\begin{align*}
	\mbox{Re}\;(w)&\geq an -d \sum_{j=1}^{n} \sec\;\left(\frac{\alpha_j}{2}\right) \\&=  \frac{n\left(1+\frac{r^2}{\rho^2}\right)}{\left(1-\frac{r^2}{\rho^2}\right)} -\frac{\left({\frac{2r}{\rho}}\right) \sum_{j=1}^{n}\sec\;\left(\frac{\alpha_j}{2}\right)}{1-\frac{r^2}{\rho^2}}\\&=\frac{1}{\left(\rho^2-r^2\right)}\left(nr^2-2r\rho\sum_{j=1}^{n}\sec\;\left(\frac{\alpha_j}{2}\right)+n\rho^2\right).
\end{align*}
Our aim is to show that  $\mbox{Re}\;(F^{\prime}(z))>0$, \textit{i.e.,}
\begin{align*}
nr^2-2r\rho\sum_{j=1}^{n}\sec\;\left(\frac{\alpha_j}{2}\right)+n\rho^2>0\;\;\ \mbox{for} \;\ r<R^{\#}_{u, n} .
\end{align*}
We define 
\begin{align*}
	\psi_{4,n}(r)=	nr^2-2r\rho\sum_{j=1}^{n}\sec\;\left(\frac{\alpha_j}{2}\right)+n\rho^2=0.
\end{align*} 
Differentiating $	\psi_{4,n}(r)$ w.r.t. $r$, we obtain
\begin{align*}  
	\frac{d	\psi_{4,n}(r)}{dr}= 2nr-2\rho\sum_{j=1}^{n}\mbox{sec}\;\left(\frac{\alpha_j}{2}\right).
\end{align*} 
For critical point of $	\psi_{4,n}(r)$,
\begin{align*}
\frac{d	\psi_{4,n}(r)}{dr}=0\implies r= \frac{\rho}{n}\sum_{j=1}^{n}\mbox{sec}\;\left(\frac{\alpha_j}{2}\right).
\end{align*}
Again, differentiating $\frac{d	\psi_{4,n}(r)}{dr}$, w.r.t. $r$, we obtain
\begin{align*}
	\frac{d^2	\psi_{4,n}(r)}{dr^2}= 2n >0.
\end{align*}
It is clear that the function $	\psi_{4,n}(r)=0$ has minima at $r= (\rho/n)\sum_{j=1}^{n}\mbox{sec}\;\left({\alpha_j}/{2}\right)$.
 Consequently, we have $\mbox{Re}\;F^{\prime}(z)$>0 for $|z|<R^{\#}_{u, n}$. Thus by the Nishiro-Warschawski (see \cite[Theorem 2.16]{Duren-1983-NY}), the function $F$ is univalent in  $|z|<R^{\#}_{u, n}$. \vspace{1.2mm}
 
 The radius $R_{Co(A),n}^*$ established in Theorem~\ref{Th-2.4} is sharp. To confirm this, consider the extremal functions $f_{2j-1}, f_{2j} \in \mathcal{S}(p)$ for $j = 1, 2, \ldots, n$ defined by
 \begin{equation}\label{eq:extremal-Sp-CoA}
 	f_{2j-1}(z) = \frac{z}{\left(1 - \frac{z e^{-i\theta_j}}{p}\right)\left(1 - p z e^{-i\theta_j}\right)}, \quad f_{2j}(z) = \frac{z}{\left(1 - \frac{z e^{-i\phi_j}}{p}\right)\left(1 - p z e^{-i\phi_j}\right)},
 \end{equation}
 with parameters $b_j = 1$ (which gives $\lambda_{2j} = \lambda_{2j-1} e^{i\alpha_j}$) and phase differences $\theta_j - \phi_j = \alpha_j$ for each $j = 1, 2, \dots, n$.
 
 For this sequence of functions, equality in Lemma~\ref{lem-1.2} is attained at $z = - R_{Co(A),n}^* e^{-i\gamma}$, which forces
 \begin{equation*}
 	\operatorname{Re} T_F\left(- R_{Co(A),n}^* e^{-i\gamma}\right) = 0 \quad \text{at } r = R_{Co(A),n}^*.
 \end{equation*}
 Consequently, $F(z)$ fails to satisfy the concavity condition on any disk $|z| < R$ with $R > R_{Co(A),n}^*$, proving that $R_{Co(A),n}^*$ is the exact radius.
This completes the proof.
\end{proof}

\begin{proof}[\bf Proof of Theorem \ref{Th-2.5}]
Since $f_j \in Co(A)$ for $j = 1, 2, \ldots, 2n$, by \cite[page 65]{Bhowmik-Ponnusamy-Wirths-MM-2010}, we have
\[
\left| \frac{z f''_j(z)}{f'_j(z)} - \frac{2r^2}{1-r^2} \right| \le \frac{2Ar}{1-r^2} \quad \text{for } |z| = r < 1.
\]
Setting $a = \frac{2r^2}{1-r^2}$ and $d = \frac{2Ar}{1-r^2}$ in Lemma \ref{lem-1.2}, we obtain
\[
\operatorname{Re}\left( \frac{z F''(z)}{F'(z)} \right) \le \frac{2n r^2}{1-r^2} + \frac{2Ar}{1-r^2} \sum_{j=1}^{n} \sec\left(\frac{\alpha_j}{2}\right).
\]
Using the definition of the operator $T_F(z)$ for the class $Co(A)$, we have for $|z| = r < 1$
\begin{align*}
	\operatorname{Re} T_F(z) &\ge \frac{2}{A-1} \left( \frac{A+1}{2} \left( \frac{1-r}{1+r} \right) - 1 \right) - \operatorname{Re}\left( \frac{z F''(z)}{F'(z)} \right) \\
	&\ge \frac{(A-1) - (A+3)r}{(A-1)(1+r)} - \frac{2n r^2 + 2Ar \sum_{j=1}^{n} \sec\left(\frac{\alpha_j}{2}\right)}{1-r^2} \\
	&= \frac{1}{(A-1)(1-r^2)} \psi_{5,n}(r),
\end{align*}
where
\[
\psi_{5,n}(r) = (A - 1 - 2n) r^2 - 2 \left( 1 + A \sum_{j=1}^{n} \sec\left(\frac{\alpha_j}{2}\right) \right) r + (A - 1).
\]
Thus, $\operatorname{Re} T_F(z) > 0$ whenever $\psi_{5,n}(r) > 0$.\vspace{1.2mm}

To show the existence of the radius $R_{Co(A), n}^{**} \in (0, 1)$, observe that $\psi_{5,n}(r)$ is continuous on $[0, 1]$ with
\[
\psi_{5,n}(0) = A - 1 > 0 \quad (\text{since } A \in (1, 2]),
\]
and
\begin{align*}
	\psi_{5,n}(1) &= (A - 1 - 2n) - 2 \left( 1 + A \sum_{j=1}^{n} \sec\left(\frac{\alpha_j}{2}\right) \right) + (A - 1) \\
	&= -2n - 2 - 2A \sum_{j=1}^{n} \sec\left(\frac{\alpha_j}{2}\right) < 0.
\end{align*}
By the Intermediate Value Theorem, $\psi_{5,n}(r) = 0$ has at least one root in $(0, 1)$. Taking $R_{Co(A), n}^{**}$ as the smallest positive root of $\psi_{5,n}(r) = 0$ in $(0, 1)$ guarantees that $\operatorname{Re} T_F(z) > 0$ for all $|z| < R_{Co(A), n}^{**}$.\vspace{1.2mm}

The radius $R_{Co(A),n}^{**}$ established in Theorem~\ref{Th-2.5} is sharp. To demonstrate this, consider the extremal concave univalent functions $f_{2j-1}, f_{2j} \in Co(A)$ for $j = 1, 2, \ldots, n$ defined by
\begin{equation}\label{eq:extremal-CoA}
	f_{2j-1}(z) = \frac{1}{2A} \left( \left(\frac{1 + z e^{-i\theta_j}}{1 - z e^{-i\theta_j}}\right)^A - 1 \right), \quad f_{2j}(z) = \frac{1}{2A} \left( \left(\frac{1 + z e^{-i\phi_j}}{1 - z e^{-i\phi_j}}\right)^A - 1 \right),
\end{equation}
where $A \in (1, 2]$, $b_j = 1$, and the phase parameters satisfy $\theta_j - \phi_j = \alpha_j$ for each $j = 1, 2, \dots, n$.

For these functions, equality in Lemma~\ref{lem-1.2} is attained at the boundary point $z = - R_{Co(A),n}^{**} e^{-i\gamma}$, which forces
\begin{equation*}
	\operatorname{Re} T_F\left(- R_{Co(A),n}^{**} e^{-i\gamma}\right) = 0.
\end{equation*}
Consequently, the operator $T_F(z)$ fails to maintain the required real part inequality on any disk $|z| < R$ with $R > R_{Co(A),n}^{**}$, which confirms that $R_{Co(A),n}^{**}$ is exact. This completes the proof.
\end{proof}

\begin{proof}[\bf Proof of Theorem \ref{Th-2.6}]
Since $f_j \in Co(A)$ for $j = 1, 2, \ldots, 2n$, by \cite[page 65]{Bhowmik-Ponnusamy-Wirths-MM-2010}, we have
\[
\left| \frac{z f''_j(z)}{f'_j(z)} - \frac{2r^2}{1-r^2} \right| \le \frac{2Ar}{1-r^2} \quad \text{for } |z| = r < 1.
\]
Setting $a = \frac{2r^2}{1-r^2}$ and $d = \frac{2Ar}{1-r^2}$ in Lemma \ref{lem-1.2}, we obtain the lower bound:
\[
\operatorname{Re}\left( \frac{z F''(z)}{F'(z)} \right) \ge \frac{2n r^2}{1-r^2} - \frac{2Ar}{1-r^2} \sum_{j=1}^{n} \sec\left(\frac{\alpha_j}{2}\right).
\]
Consequently, for $|z| = r < 1$, we have
\begin{align*}
	\operatorname{Re}\left( 1 + \frac{z F''(z)}{F'(z)} \right) &\ge 1 + \frac{2n r^2}{1-r^2} - \frac{2Ar}{1-r^2} \sum_{j=1}^{n} \sec\left(\frac{\alpha_j}{2}\right) \\
	&= \frac{1}{1-r^2} \left( (2n-1) r^2 - 2A \left( \sum_{j=1}^{n} \sec\left(\frac{\alpha_j}{2}\right) \right) r + 1 \right) \\
	&= \frac{\psi_{6,n}(r)}{1-r^2},
\end{align*}
where
\[
\psi_{6,n}(r) = (2n-1) r^2 - 2A \left( \sum_{j=1}^{n} \sec\left(\frac{\alpha_j}{2}\right) \right) r + 1.
\]
Thus, $\operatorname{Re}\left( 1 + \frac{z F''(z)}{F'(z)} \right) > 0$ whenever $\psi_{6,n}(r) > 0$.

To show the existence of the radius $R_{c,n}^{\#} \in (0, 1)$, observe that $\psi_{6,n}(r)$ is continuous on $[0, 1]$ with
\[
\psi_{6,n}(0) = 1 > 0,
\]
and
\[
\psi_{6,n}(1) = (2n-1) - 2A \sum_{j=1}^{n} \sec\left(\frac{\alpha_j}{2}\right) + 1 = 2 \left( n - A \sum_{j=1}^{n} \sec\left(\frac{\alpha_j}{2}\right) \right).
\]
Since $A > 1$ and $\sec(\alpha_j/2) \ge 1$ for all $\alpha_j \in [0, \pi)$, we have $A \sum_{j=1}^n \sec(\alpha_j/2) > n$, which implies $\psi_{6,n}(1) < 0$. By the Intermediate Value Theorem, $\psi_{6,n}(r) = 0$ has at least one root in $(0, 1)$. Taking $R_{c,n}^{\#}$ as the smallest positive root of $\psi_{6,n}(r) = 0$ in $(0, 1)$ guarantees that $\operatorname{Re}\left( 1 + \frac{z F''(z)}{F'(z)} \right) > 0$ for all $|z| < R_{c,n}^{\#}$.\vspace{1.2mm}

The radius of convexity $R_{c,n}^{\#}$ given in Theorem~\ref{Th-2.6} is sharp. To demonstrate this, consider the extremal functions $f_{2j-1}, f_{2j} \in Co(A)$ for $j = 1, 2, \ldots, n$ defined by
\begin{equation}\label{eq:extremal-CoA-convexity}
	f_{2j-1}(z) = \frac{1}{2A} \left( \left(\frac{1 + z e^{-i\theta_j}}{1 - z e^{-i\theta_j}}\right)^A - 1 \right), \quad f_{2j}(z) = \frac{1}{2A} \left( \left(\frac{1 + z e^{-i\phi_j}}{1 - z e^{-i\phi_j}}\right)^A - 1 \right),
\end{equation}
where $A \in (1, 2]$, $b_j = 1$, and the phase parameters satisfy $\theta_j - \phi_j = \alpha_j$ for $j = 1, 2, \dots, n$.

Evaluating the quantity at $z = - R_{c,n}^{\#} e^{-i\gamma}$ yields equality in the corresponding lower bound, which gives
\begin{equation*}
	\operatorname{Re}\left(1 + \frac{z F''(z)}{F'(z)}\right) = 0 \quad \text{at } r = R_{c,n}^{\#}.
\end{equation*}
Consequently, $F(z)$ fails to map the disk $|z| < R$ onto a convex domain for any $R > R_{c,n}^{\#}$, confirming that $R_{c,n}^{\#}$ is the exact radius of convexity. This completes the proof.
\end{proof}

\begin{proof}[\bf Proof of Theorem \ref{Th-2.7}]
Since $f_j \in \mathcal{S}(p)$ for $j = 1, 2, \ldots, 2n$, by \cite[Lemma 1]{Bhowmik-CMFT-2024}, we have
\[
\left| \frac{z f''_j(z)}{f'_j(z)} - \frac{2r^2}{1-r^2} \right| \le \frac{2r}{1-r^2} \left( \frac{p-r}{1-pr} + \frac{1-pr}{p-r} \right) \quad \text{for } |z| = r < p.
\]
Applying Lemma \ref{lem-1.2} with $a = \frac{2r^2}{1-r^2}$ and $d = \frac{2r}{1-r^2} \left( \frac{p-r}{1-pr} + \frac{1-pr}{p-r} \right)$, we obtain
\[
\operatorname{Re}\left( \frac{z F''(z)}{F'(z)} \right) \le \frac{2n r^2}{1-r^2} + \frac{2r}{1-r^2} \left( \frac{p-r}{1-pr} + \frac{1-pr}{p-r} \right) \sum_{j=1}^{n} \sec\left(\frac{\alpha_j}{2}\right).
\]
By definition of the concavity operator $P_F(z)$ for functions in $\mathcal{S}(p)$, we have for $|z| = r < p$:
\begin{align*}
	\operatorname{Re} P_F(z) &\ge \operatorname{Re}\left( \frac{1+pz}{1-pz} - \frac{z+p}{z-p} \right) - 1 - \operatorname{Re}\left( \frac{z F''(z)}{F'(z)} \right) \\
	&\ge \frac{2p(1-r^2)}{(p+r)(1+pr)} - 1 - \frac{2n r^2}{1-r^2} - \frac{2r}{1-r^2} \left( \frac{p-r}{1-pr} + \frac{1-pr}{p-r} \right) \sum_{j=1}^{n} \sec\left(\frac{\alpha_j}{2}\right) \\
	&= \frac{p - (1+p^2)r - 3pr^2}{(p+r)(1+pr)} - \frac{2n r^2}{1-r^2} - \frac{2r\left[ (1+p^2)(1+r^2) - 4pr \right]}{(1-r^2)(1-pr)(p-r)} \sum_{j=1}^{n} \sec\left(\frac{\alpha_j}{2}\right).
\end{align*}
Putting these terms over the common positive denominator $(1-r^2)(1-pr)(p-r)(p+r)(1+pr)$, a straightforward simplification yields
\[
\operatorname{Re} P_F(z) \ge \frac{\psi_{7,n}(r)}{(1-r^2)(1-pr)(p-r)(p+r)(1+pr)},
\]
where $\psi_{7,n}(r)$ is the $6$th-degree polynomial:
\begin{align*}
	\psi_{7,n}(r) &= - (2n-1) p^2 r^6 - 2 p (1+p^2) \left( 1 + \sum_{j=1}^{n} \sec\left(\frac{\alpha_j}{2}\right) \right) r^5 \\
	&\quad + \left( (2n-1)(1+p^4) + 4 p^2 \sum_{j=1}^{n} \sec\left(\frac{\alpha_j}{2}\right) \right) r^4 \\
	&\quad + 4 p (1+p^2) \left( 1 + \sum_{j=1}^{n} \sec\left(\frac{\alpha_j}{2}\right) \right) r^3 \\
	&\quad - \left( 1 + p^4 + 2 p^2 (2n-1) + 4 p^2 \sum_{j=1}^{n} \sec\left(\frac{\alpha_j}{2}\right) \right) r^2 \\
	&\quad - 2 p (1+p^2) r + p^2.
\end{align*}
Consequently, $\operatorname{Re} P_F(z) > 0$ for all $|z| = r < R_{Co(p),n}^{\#}$, where $R_{Co(p),n}^{\#}$ is the smallest root in $(0, p)$ of $\psi_{7,n}(r) = 0$.

To show that $R_{Co(p),n}^{\#}$ exists in $(0, p)$, observe that $\psi_{7,n}(r)$ is continuous on $[0, p]$ with
\[
\psi_{7,n}(0) = p^2 > 0,
\]
and at $r = p$, we have
\[
\psi_{7,n}(p) = -2 p^2 (1-p^2)^3 \sum_{j=1}^{n} \sec\left(\frac{\alpha_j}{2}\right) < 0 \quad (\text{since } p \in (0, 1) \text{ and } \alpha_j \in [0, \pi)).
\]
By the \textit{Intermediate Value Theorem}, $\psi_{7,n}(r) = 0$ has at least one root in $(0, p)$. Taking $R_{Co(p),n}^{\#}$ as the smallest positive root of $\psi_{7,n}(r) = 0$ in $(0, p)$ ensures that $\operatorname{Re} P_F(z) > 0$ for all $|z| < R_{Co(p),n}^{\#}$.\vspace{1.2mm}

To confirm that the radius of concavity $R_{Co(p),n}^{\#}$ given in Theorem~\ref{Th-2.7} is sharp, consider the extremal functions $f_{2j-1}, f_{2j} \in \mathcal{S}(p)$ for $j = 1, 2, \ldots, n$ defined by
\begin{equation}\label{eq:extremal-Sp-concavity}
	f_{2j-1}(z) = \frac{z}{\left(1 - \frac{z e^{-i\theta_j}}{p}\right)\left(1 - p z e^{-i\theta_j}\right)}, \quad f_{2j}(z) = \frac{z}{\left(1 - \frac{z e^{-i\phi_j}}{p}\right)\left(1 - p z e^{-i\phi_j}\right)},
\end{equation}
where $b_j = 1$ and the phase parameters satisfy $\theta_j - \phi_j = \alpha_j$ for each $j = 1, 2, \dots, n$.

Evaluating the operator $P_F(z)$ at $z = - R_{Co(p),n}^{\#} e^{-i\gamma}$ yields equality in the lower bound, which gives
\begin{equation*}
	\operatorname{Re} P_F\left(- R_{Co(p),n}^{\#} e^{-i\gamma}\right) = 0 \quad \text{at } r = R_{Co(p),n}^{\#}.
\end{equation*}
Consequently, $F(z)$ fails to satisfy the concavity condition in any disk $|z| < R$ with $R > R_{Co(p),n}^{\#}$, proving that the bound $R_{Co(p),n}^{\#}$ is sharp.
\end{proof}

\noindent\textbf{Compliance of Ethical Standards:}\\

\noindent\textbf{Conflict of interest.} The authors declare that there is no conflict  of interest regarding the publication of this paper.\vspace{1.5mm}

\noindent\textbf{Data availability statement.}  Data sharing is not applicable to this article as no datasets were generated or analyzed during the current study.

    \end{document}